\documentclass[12pt]{amsart}

\usepackage{amsmath}
\usepackage{amssymb}
\usepackage{array}

\usepackage[T2A]{fontenc}

\newtheorem{theorem}{Теорема}
\newtheorem{lemma}{Лемма}
\newtheorem{corollary}{Следствие}

\newtheorem{conjecture}{Предположение}

\begin{document}


\vspace{0.3cm}

\title{Рациональные приближения значений дигамма-функции  и
гипотеза о знаменателях}

\author{Т.~Хессами Пилеруд, Х.~Хессами Пилеруд }

\address{Mathematics Department, Faculty of Basic
Sciences,  Shahrekord University, \-\-Shahrekord, P.O.~Box~115,    Iran.}
\email{hessamik@gmail.com, hessamit@ipm.ir,
hessamit@gmail.com}





\begin{abstract}
In 2007, A.~I.~Aptekarev and his collaborators discovered a sequence of rational approximations
to Euler's constant $\gamma$ defined by a  linear recurrence.
In this paper, we generalize this result and
present an  explicit construction of rational approximations for the numbers
$\ln(b)-\psi(a+1),$ $a, b\in {\mathbb Q},$ $b>0, a>-1,$ where $\psi$ defines the logarithmic
derivative of the Euler gamma function.
We prove exact formulas for denominators and numerators of the approximations in terms of hypergeometric sums.
As a consequence, we get rational approximations for the numbers $\pi/2\pm\gamma.$ We compare the results obtained
with those  of T.~Rivoal for the numbers $\gamma+\ln(b)$ and prove denominators conjectures
proposed by Rivoal for denominators of rational approximations for $\gamma+\ln(b)$ and common denominators of simultaneous approximations for the numbers $\gamma$ and $\zeta(2)-\gamma^2.$

\end{abstract}

\maketitle

\section{Введение}

{\bf 1.} Пусть $a, b\in {\mathbb R},$ $b>0, a>-1.$ Рассмотрим многочлен Якоби-Лагерра
\begin{equation}
Q_n^{(a,b)}(x)=\frac{1}{n!^2}\frac{x^{-a}e^{bx}}{1-x}
(x^n(x^{n+a}(1-x)^{2n+1}e^{-bx})^{(n)})^{(n)},
\label{eq01}
\end{equation}
являющийся совместно ортогональным многочленом относительно двух весовых функций
$$
w_1(x)=x^a(1-x)e^{-bx}, \qquad w_2(x)=x^a(1-x)e^{-bx}\ln(x)
$$
на промежутках $\triangle_1=[0,1]$  и $\triangle_2=[1,+\infty),$ и для которого имеют
место следующие соотношения ортогональности
$$
\int_{\triangle_j}x^{\nu}Q_n^{(a,b)}(x)w_k(x)\,dx=0, \qquad \nu=0,1,\ldots, n-1, \quad j,k=1,2.
$$
Совместно ортогональные многочлены (\ref{eq01}) были введены в статье \cite{as}.
В работе \cite{ap} многочлены $Q_n^{(a,b)}$ при  $a=0, b=1$ впервые были применены для построения
рациональных приближений постоянной Эйлера
$$
\gamma=\lim_{k\to\infty}\left(1+\frac{1}{2}+\cdots+\frac{1}{k}-\ln k\right).
$$
Полученные приближения хотя и не позволяют доказать иррациональность $\gamma,$ об\-ла\-да\-ют
достаточно хорошими свойствами. Числители $p_n$ и знаменатели   $q_n$ ра\-цио\-наль\-ных приближений
являются целыми числами, задаваемыми рекурсией третьего порядка
\begin{equation}
\begin{split}
(16n-15)q_{n+1}=&(128n^3+40n^2-82n-45)q_n-n^2(256n^3-240n^2+64n-7)q_{n-1} \\
&+n^2(n-1)^2(16n+1)q_{n-2}
\label{eq02}
\end{split}
\end{equation}
с начальными условиями
\begin{equation*}
\begin{array}{ccc}
p_0=0, \qquad & \qquad p_1=2, \qquad & \qquad p_2=31, \\
q_0=1,  \qquad &  \qquad q_1=3, \qquad &  \qquad q_2=50,
\end{array}
\end{equation*}
и имеющими следующую асимптотику
\begin{equation}
\begin{split}
q_n&=(2n)!\frac{e^{\sqrt{2n}}}{\sqrt[4]{n}}\left(\frac{1}{\sqrt{\pi}(4e)^{3/8}}+O(n^{-1/2})\right),\\[3pt]
p_n-\gamma q_n&=(2n)!\frac{e^{-\sqrt{2n}}}{\sqrt[4]{n}}\left(\frac{2\sqrt{\pi}}{(4e)^{3/8}}+O(n^{-1/2})\right).
\label{eq03}
\end{split}
\end{equation}
Таким образом, рекуррентные соотношения (\ref{eq02})  определяют рациональные при\-бли\-же\-ния постоянной Эйлера
$$
\frac{p_n}{q_n}-\gamma=2\pi e^{-2\sqrt{2n}}\left(1+O(n^{-1/2})\right), \qquad n\to\infty.
$$
Целочисленность последовательностей $p_n$ и $q_n$  не следует непосредственно из (\ref{eq02})
и была доказана независимо в \cite{tu}
с помощью рассмотрения более ``плотной'' по\-сле\-до\-ва\-тель\-нос\-ти рациональных приближений к $\gamma,$
числители и знаменатели которой удовлетворяют более простым рекуррентным соотношениям.
Основную трудность работ \cite{ap, tu} представляет нахождение рекуррентных уравнений, которые
опираются на ре\-кур\-рент\-ные соотношения для многочленов (\ref{eq01}).

В настоящей статье мы обобщим конструкцию работы \cite{ap} и
 построим ра\-цио\-наль\-ные приближения для  чисел $\ln b-\psi(a+1),$ где $a,b\in {\mathbb Q},
a>-1, b>0,$  $\psi(z)$ определяет логарифмическую производную гамма-функции Эйлера (или дигамма-функцию)
$$
\psi(z)=\frac{d}{dz}\log\Gamma(z)=\frac{\Gamma'(z)}{\Gamma(z)}=-\gamma+\sum_{n=0}^{\infty}
\left(\frac{1}{n+1}-\frac{1}{n+z}\right),
$$
при этом $\psi(1)=-\gamma.$ Мы дадим явные формулы для числителей и знаменателей этих приближений
в терминах гипергеометрических сумм,  а для остатков приближений --- в терминах
гиперэкспоненциальных интегралов. Хорошо развитая техника для работы с такими объектами
(см. \cite{ab, maz}), реализованная в алгоритмах EKHAD и MultiAlmkvistZeilberger, сопровождающих
программный пакет MAPLE, позволяет получать рекуррентные уравнения, задающие рациональные приближения,
ав\-то\-ма\-ти\-чес\-ки.

Сформулируем результаты, доказываемые в статье. Пусть, как обычно, $(a)_0=1,$ $(a)_n=a(a+1)\ldots (a+n-1),$
$n\ge 1,$ обозначает символ Похгаммера или сдвинутый факториал.
Гипергеометрическая функция Гаусса $F(a,b;c;z)$ при комплексных параметрах $a, b, c$ определяется
как сумма ряда
\begin{equation}
F(a,b;c;z)=\sum_{k=0}^{\infty}\frac{(a)_k(b)_k}{k!(c)_k}z^k, \qquad
c\ne 0, -1, -2, \ldots.
\label{eqf}
\end{equation}
Если $a$ и $b$ отличны от $0,-1,-2,\ldots,$ то гипергеометрический ряд (\ref{eqf}) абсолютно сходится
для всех значений $z,$ лежащих внутри круга $|z|<1.$ Если один из па\-ра\-мет\-ров $a, b$
равен целому неположительному числу, то ряд (\ref{eqf}) обрывается на конечном числе слагаемых,
и гипергеометрическая функция является многочленом от $z.$  Пусть $H_n(\alpha)$ обозначает  частичную сумму
гармонического ряда
$$
H_n(\alpha)=\sum_{k=1}^n\frac{1}{k+\alpha}, \quad \alpha\ne -1, \ldots, -n, \,\, n\ge 1, \quad
H_0(\alpha)=0,
$$
$H_n:=H_n(0),$ и
 $D_n$ --- наименьшее общее кратное чисел
$1,2, \ldots, n.$ Всюду в даль\-ней\-шем считаем, что сумма (произведение) по пустому множеству
значений индекса равна 0 (равно $1$).

\begin{theorem}\label{t1}
Для любого целого $n\ge 0$ справедливо следующее представление:
$$
\frac{b^{2n+a+1}}{\Gamma(a+1)}\int_0^{\infty}Q_n^{(a,b)}(x)x^ae^{-bx}\ln(x)\,dx=
p_n-q_n(\ln(b)-\psi(a+1)),
$$
где
$$
q_n=\sum_{k=0}^n\binom{n}{k}^2(a+1)_{n+k}\,b^{n-k}\in {\mathbb Z}[a,b],
$$
$$
p_n=\sum_{k=0}^n\binom{n}{k}^2(a+1)_{n+k}\,b^{n-k}H_{n+k}(a) 
$$
$$
-2\sum_{k=1}^n\sum_{m=0}^k\sum_{l=0}^{n-k}\frac{(-1)^{m+k}}{k}\binom{n}{k+l}\binom{k}{m}
\binom{n-k}{l}(a+1)_{m+n+l}\,b^{n-m-l}
\in \frac{1}{D_n}{\mathbb Z}[a,b].
$$
При этом для последовательностей $q_n$ и $p_n$ имеют место следующие асимп\-то\-ти\-чес\-кие формулы:
$$
q_n=(2n)!\,\frac{e^{\sqrt{2bn}}}{n^{1/4-a}}(c(a,b)+O(n^{-1/2})), \quad n\to\infty
$$
$$
p_n-q_n(\ln b-\psi(a+1))=(2n)!\,\frac{-e^{\sqrt{2bn}}}{n^{1/4-a}}(2\pi c(a,b)+O(n^{-1/2})),
\quad n\to\infty
$$
где
 $$
 c(a,b)=2^a/\Gamma(a+1)/(8b\pi^2e^{3b/2})^{1/4}.
 $$
\end{theorem}
Последовательности $q_n$ и $p_n$  из теоремы \ref{t1} являются решениями рекуррентного уравнения третьего
порядка с полиномиальными коэффициентами, зависящими от $a$ и $b.$ Это рекуррентное уравнение легко
может быть выписано при конкретных значениях $a$ и $b$ (см.~замечание к лемме \ref{l4}).
Однако, мы не приводим его здесь в общем виде в силу громоздкости его коэффициентов.

Заметим, что при $a=0, b=1$ полученные приближения совпадают с ра\-цио\-наль\-ными приближениями Аптекарева
 для константы Эйлера. Тем самым, теорема 1 дает явные формулы для $q_n, p_n$
из (\ref{eq03})
\begin{equation}
q_n=\sum_{k=0}^n\binom{n}{k}^2(n+k)!,
\label{znam}
\end{equation}
$$
p_n=\sum_{k=0}^n\binom{n}{k}^2(n+k)!H_{n+k}
-2\sum_{k=1}^n\sum_{m=0}^k\sum_{l=0}^{n-k}\frac{(-1)^{m+k}}{k}\binom{n}{k+l}\binom{k}{m}\binom{n-k}{l}(m+n+l)!.
$$
Откуда следует не только целочисленность последовательностей $q_n$ и $p_n,$ но также и делимость $q_n$ на $n!,$
а $p_n$  на $\frac{n!}{D_n}.$ Поэтому линейные формы (\ref{eq03}) можно сократить на большой общий множитель
$\frac{n!}{D_n}:$ полученные формы  также будут иметь целые коэффициенты, однако, этого все равно недостаточно
для доказательства ир\-ра\-цио\-наль\-нос\-ти $\gamma.$

Отметим несколько следствий теоремы \ref{t1}. При $a=0,$ $b\in {\mathbb Q}, b>0$ получаем
рациональные приближения для $\gamma+\ln(b).$
\begin{corollary} \label{c1}
Последовательности
$$
q_n=\sum_{k=0}^n\binom{n}{k}^2(n+k)!\,b^{n-k}\in n!\,{\mathbb Z}[b]
$$
и
$$
p_n=\sum_{k=0}^n\binom{n}{k}^2(n+k)!\,b^{n-k}H_{n+k}
$$
$$
-2\sum_{k=1}^n
\sum_{m=0}^k\sum_{l=0}^{n-k}\frac{(-1)^{m+k}}{k}\binom{n}{k+l}\binom{k}{m}\binom{n-k}{l}
(m+n+l)!\,b^{n-m-l}\in \frac{n!}{D_n}\,{\mathbb Z}[b]
$$
удовлетворяют следующему асимптотическому равенству
$$
\frac{p_n}{q_n}-(\ln(b)+\gamma)=2\pi e^{-2\sqrt{2bn}}(1+O(n^{-1/2})), \qquad n\to\infty.
$$
\end{corollary}
{\it Замечание.} При $b\not\in \frac{4}{3}+\frac{16}{3}{\mathbb N}$
последовательности $p_n$ и $q_n$ из следствия \ref{c1}
 аль\-тер\-на\-тив\-но могут быть заданы как
решения рекуррентного уравнения третьего порядка
\begin{equation*} \begin{split} &(16n-3b+4)f_{n+2} \\ &=(128n^3+8n^2b-6nb^2+416n^2-12nb-9b^2+400n-38b+88)f_{n+1}\\ &-(n+1)^2(256n^3-48n^2b+16nb^2-3b^3+576n^2-48nb+8b^2+384n-12b+80)f_{n}\\ &+b^2n^2(n+1)^2(16n-3b+20)f_{n-1} \end{split} \end{equation*}
с начальными условиями
\begin{equation*}
\begin{array}{ccc}
p_0=0, \qquad & \qquad p_1=3b-1, \qquad & \qquad p_2=9b^2+44b-22, \\
q_0=1,  \qquad &  \qquad q_1=b+2, \qquad &  \qquad q_2=2b^2+24b+24.
\end{array}
\end{equation*}
При $b=1,$ $a\in {\mathbb Q},$ $a>-1$ получаем рациональные приближения зна\-че\-ний дигамма-функции
$\psi(a+1).$
\begin{corollary} \label{c2}
При $a\in {\mathbb Q},$ $a>-1$ имеем
$$
\psi(a+1)-\frac{p_n}{q_n}=2\pi e^{-2\sqrt{2n}}(1+O(n^{-1/2})), \qquad n\to\infty
$$                                                                                                                                                                                                                                                                                                            где
$$   q_n=\sum_{k=0}^n\binom{n}{k}^2(a+1)_{n+k},     $$
\begin{equation*}
\begin{split}
p_n=2\sum_{k=1}^n\sum_{m=0}^k\sum_{l=0}^{n-k}&\frac{(-1)^{m+k}}{k}\binom{n}{k+l}  \binom{k}{m}\binom{n-k}{l}(a+1)_{m+n+l} \\ &-\sum_{k=0}^n\binom{n}{k}^2(a+1)_{n+k}
H_{n+k}(a).
\end{split}
\end{equation*}
\end{corollary}

При $b=1/8,$ $a=-3/4,$ $a=-1/4$ получаем рациональные приближения для чисел $\pi/2\pm\gamma.$
\begin{corollary} \label{c3}
Последовательности
$$
q_n=\sum_{k=0}^n\binom{n}{k}^2\left(\frac{1}{4}\right)_{n+k}8^{n+k},
$$
\begin{equation*}
\begin{split}
p_n&=\sum_{k=0}^n\binom{n}{k}^2\left(\frac{1}{4}\right)_{n+k}\sum_{l=1}^{n+k}
\frac{4\cdot 8^{n+k}}{4l-3}\\
&-2\sum_{k=1}^n
\sum_{m=0}^k\sum_{l=0}^{n-k}\frac{(-1)^{m+k}}{k}\binom{n}{k+l}\binom{k}{m}\binom{n-k}{l}
\left(\frac{1}{4}\right)_{m+n+l}8^{m+n+l}
\end{split}
\end{equation*}
задают рациональные приближения числа $\gamma+\pi/2$
$$
\frac{p_n}{q_n}-\left(\gamma+\frac{\pi}{2}\right)=2\pi e^{-\sqrt{n}}(1+O(n^{-1/2})), \quad n\to\infty.
$$
\end{corollary}
{\it Замечание.}  Последовательности $q_n$ и $p_n$ из следствия \ref{c3}
также могут быть заданы как решения
рекуррентного уравнения
\begin{equation*}
\begin{split}
&(2n-1)(n+2)(2048n^3-304n^2-456n+153)f_{n+2}\\
&=(2097152n^7+7618560n^6+7195136n^5
-1843584n^4-3710560n^3+569496n^2\\
&+505524n-142182)f_{n+1}-
4(n+1)(67108864n^8+191365120n^7+149225472n^6\\
&-2504192n^5-43117824n^4-10174400n^3+2244800n^2+726450n+45729)f_n\\
&+256(2n+1)(n+1)(2048n^3+5840n^2+5080n+1441)(4n-3)^2n^2f_{n-1}
\end{split}
\end{equation*}
с начальными условиями 
$$
p_1=0, \quad p_1=68, \quad p_2=31596\quad\text{и}\quad q_0=1, \quad q_1=22, \quad q_2=10820.
$$
\begin{corollary} \label{c4}
Последовательности
$$
q_n=\sum_{k=0}^n\binom{n}{k}^2\left(\frac{3}{4}\right)_{n+k}8^{n+k},
$$
\begin{equation*}
\begin{split}
p_n&=2\sum_{k=1}^n
\sum_{m=0}^k\sum_{l=0}^{n-k}\frac{(-1)^{m+k}}{k}\binom{n}{k+l}\binom{k}{m}\binom{n-k}{l}
\left(\frac{3}{4}\right)_{m+n+l}8^{m+n+l}\\
&-\sum_{k=0}^n\binom{n}{k}^2\left(\frac{3}{4}\right)_{n+k}\sum_{l=1}^{n+k}
\frac{4\cdot 8^{n+k}}{4l-1}
\end{split}
\end{equation*}
задают рациональные приближения числа $\pi/2-\gamma$
$$
\left(\frac{\pi}{2}-\gamma\right)-\frac{p_n}{q_n}=2\pi e^{-\sqrt{n}}(1+O(n^{-1/2})), \quad n\to\infty.
$$
\end{corollary}
{\it Замечание.}  Последовательности $q_n$ и $p_n$ из следствия \ref{c4}
являются решениями
рекуррентного уравнения
\begin{equation*}
\begin{split}
&(n+2)(4096n^4+2464n^3-736n^2+58n-3)f_{n+2}\\
&=(2097152n^7+11288576n^6+21899776n^5+17372032n^4+3526048n^3-1084856n^2\\
&-18076n+4298)f_{n+1}
-4(n+1)(67108864n^8+342360064n^7+664600576n^6\\
&+613444096n^5+259912960n^4+27448384n^3
-10104768n^2-1517358n+191079)f_n\\
&+256(n+1)(4096n^4+18848n^3+31232n^2+22362n+5879)(4n-1)^2n^2f_{n-1}
\end{split}
\end{equation*}
с начальными условиями 
$$
p_1=0, \quad p_1=-12, \quad p_2=8596\quad\text{и}\quad q_0=1, \quad q_1=90, \quad q_2=62916.
$$
{\bf 2.} В работе \cite{ri} Т.~Ривоаль предложил альтернативный метод для получения ра\-цио\-наль\-ных приближений
значений $\gamma+\ln(x),$ $x\in {\mathbb Q}, x>0.$ Для этой цели Ривоаль использовал обобщенные многочлены
Лагерра
$$
A_n(x)=\frac{1}{n!^2}e^x(x^n(x^ne^{-x})^{(n)})^{(n)},
$$
являющиеся совместно ортогональными многочленами на $[0,+\infty)$ относительно ве\-со\-вых функций $w_1(x)=e^{-x}$
и $w_2(x)=e^{-x}\ln(x).$ Построенные им приближения имеют вид
\begin{equation}
 (n+1)^2\int_0^{\infty}\frac{A_{n+1}(x)A_n(t)-A_n(x)A_{n+1}(t)}{t-x}(\ln(t)-\ln(x))e^{-t}\, dt=
Q_n(x)(\ln(x)+\gamma)-P_n(x),
\label{pq}
\end{equation}

\vspace{0.4cm}

\noindent где $n!^2P_n(x), n!^2Q_n(x)\in {\mathbb Z}[x],$ $\deg P_n, Q_n\le n+1,$
$Q_n(x)=O(\exp(3x^{1/3}n^{2/3}-x^{2/3}n^{1/3})),$ и
$$
\left|\ln(x)+\gamma-\frac{P_n(x)}{Q_n(x)}\right|=O(\exp(-9/2x^{1/3}n^{2/3}+3/2x^{2/3}n^{1/3})).
$$
В \cite{ri} Т.~Ривоаль высказал предположение, что для знаменателей $Q_n(x)$ имеет место более строгое
включение.
\begin{conjecture}
Для любого целого $n\ge 0,$ $n!\,Q_n(x)\in {\mathbb Z}[x].$
\end{conjecture}
Аналогичное предположение было высказано в \cite{ri} для по\-сле\-до\-ва\-тель\-нос\-ти $b_n,$
представляющей общие знаменатели совместных приближений чисел $\gamma$ и $\zeta(2)-\gamma^2$
по формулам
\begin{align}
(n+1)^3\!\!\!\int_0^{\infty}
\frac{\tilde{A}_n(t)\tilde{A}_{n+1}(1)-\tilde{A}_{n+1}(t)\tilde{A}_n(1)}{1-t}
\ln(t) e^{-t}\,dt&=a_{1,n}-b_n\gamma, \label{eqri1}\\[5pt]
(n+1)^3\!\!\!\int_0^{\infty} \frac{\tilde{A}_n(t)\tilde{A}_{n+1}(1)-\tilde{A}_{n+1}(t)\tilde{A}_n(1)}{1-t}
(\ln^2(t)+2\gamma\ln(t)) e^{-t}\,dt&=b_n(\zeta(2)-\gamma^2)-a_{2,n},
\label{eqri2}
\end{align}
где
$$
\tilde{A}_n(x)=\frac{1}{n!^3}e^x(x^n(x^n(x^ne^{-x})^{(n)})^{(n)})^{(n)} $$
и  $a_{1,n}, a_{2,n}, b_n\in \frac{1}{(3n)!(3n+2)!}{\mathbb Z}.$
\begin{conjecture}
Для любого целого $n\ge 0,$ $n!^2\,b_n\in {\mathbb Z}.$
\end{conjecture}
В настоящей статье мы доказываем более общее утверждение, из которого сле\-ду\-ют оба предположения Ривоаля,
а также более точные оценки  для знаменателей рас\-смат\-ри\-вае\-мых последовательностей.
\begin{theorem} \label{t2}
Пусть $s\in {\mathbb N}$  и $(A_{s,n})_{n\ge 0}$ последовательность обобщенных многочленов Лагерра,
заданная рекурсивно:
$$
A_{0,n}(x)\equiv 1, \qquad A_{r+1,n}(x)=\frac{1}{n!}e^x(x^ne^{-x}A_{r,n}(x))^{(n)}, \quad r\ge 0.
$$
Тогда для интегралов
$$
J_{s,n}(x):=(n+1)^s\!\!\int_0^{\infty}\frac{A_{s,n+1}(x)A_{s,n}(t)-A_{s,n}(x)A_{s,n+1}(t)}{t-x}\ln(t)e^{-t}\,dt
$$
справедливо следующее представление
$$
J_{s,n}(x)=Q_{s,n}(x)\gamma-P_{s,n}(x),
$$
где
$$
Q_{s,n}(x)=(n+1)^s\!\!\int_0^{\infty}\frac{A_{s,n+1}(t)A_{s,n}(x)-A_{s,n}(t)A_{s,n+1}(x)}{t-x}e^{-t}\,dt
$$

\noindent и  $n!^{s-1}Q_{s,n}(x), D_n n!^{s-1} P_{s,n}(x)\in {\mathbb Z}[x].$
\end{theorem}
\begin{corollary} \label{c5}
Для многочленов $Q_n(x), P_n(x),$ определенных в {\rm (\ref{pq}),}  имеют место следующие включения:
$n!\,Q_n(x), \,\, D_n n!P_n(x)\in{\mathbb Z}[x].$
 \end{corollary}
\begin{corollary} \label{c6}
Для последовательностей $b_n, a_{1,n}, a_{2,n},$ определенных в {\rm  (\ref{eqri1}), (\ref{eqri2})}  имеем
$n!^2 b_n, \,\, D_n n!^2 a_{1,n}, \,\, D_n^2 n!^2 a_{2,n}\in {\mathbb Z}.$
\end{corollary}

Аналогичный подход применим также к работе Ривоаля \cite{ri1}, где с помощью обобщенных
многочленов  Лагерра вида
$$
\hat{A}_{n,a}(x)=\frac{1}{n!^2}e^x(x^{n-a}(x^{n+a}e^{-x})^{(n)})^{(n)}
$$
были получены рациональные приближения для числа $\frac{\Gamma(a+1)}{x^a},$ $a,x\in {\mathbb Q},$
$a>-1,$ $x>0.$ В этой конструкции числители и знаменатели рациональных приближений
определяются из разложения следующего интеграла
\begin{equation}
\int_0^{\infty}(\hat{A}_{n,a}(x)\hat{A}_{n+1,a}(t)-\hat{A}_{n+1,a}(x)\hat{A}_{n,a}(t))
\frac{x^a-t^a}{x-t}e^{-t}\,dt=Q_n(a,x)\Gamma(1+a)-x^aP_n(a,x),
\label{eq100}
\end{equation}
при этом
$$
\left|\frac{\Gamma(a+1)}{x^a}-\frac{P_n(a,x)}{Q_n(a,x)}\right|=O(e^{-9/2x^{1/3}n^{2/3}+3/2x^{2/3}n^{1/3}}),
\quad n\to\infty,
$$
и
$$
n!^2(n+1)!^2({\rm den}\, a)^{2n+1}({\rm den}\, x)^{4n-1}\!P_n(a,x),
\,
n!^2(n+1)!^2({\rm den}\, a)^{3n+2}({\rm den}\, x)^{4n-1}Q_n(a,x)\in {\mathbb Z},
$$
где ${\rm den}\, \alpha\in {\mathbb N}$ обозначает знаменатель несократимой дроби $\alpha.$
В данной ситуации также возможно существенно уточнить оценки знаменателей
$Q_n(a,x)$ и $P_n(a,x).$

\noindent {\bf Предложение 1.} {\it При  $a,x\in {\mathbb Q}, a>-1, x>0$  справедливы следующие включения:
$$
({\rm den}\, x)^{n+1}\mu_{a}^{2n+1}P_n(a,x), \quad
({\rm den}\, x)^{n+1}\mu_{a}^{3n+1}Q_n(a,x)\in
\frac{1}{(n+1)^2 n!}
{\mathbb Z,}
$$
где
$\mu^n_{a}=({\rm den}\, a)^n\cdot\prod_{p|{\rm den}\, a}p^{[\frac{n}{p-1}]}.$
}

\vspace{0.3cm}

\section{Вспомогательные леммы.}

\begin{lemma} \label{l1}
Пусть $a, b\in {\mathbb R},$ $b>0,$ $a>-1.$ Справедливы следующие равенства:
\begin{equation*}
\begin{array}{lcl}
(i)& \qquad\qquad\qquad\qquad \displaystyle\int_0^{\infty}x^ae^{-bx}\,dx&=\displaystyle\frac{\Gamma(a+1)}{b^{a+1}},\\[15pt]
(ii)&  \qquad\qquad\qquad\qquad \displaystyle\int_0^{\infty}x^ae^{-bx}\ln(x)\,dx
&=\displaystyle\frac{\Gamma(a+1)}{b^{a+1}}(\psi(a+1)-\ln b),\\[15pt]
(iii)&  \qquad\qquad\qquad\qquad \displaystyle\int_0^{\infty}x^ae^{-x}\ln(x)^2\,dx
&=\displaystyle\Gamma(a+1)(\psi^2(a+1)+\psi'(a+1))).
\end{array}
\end{equation*}
\end{lemma}
{\bf Доказательство.} Равенство $(i)$ непосредственно следует из определения гамма-функции Эйлера.
Дифференцируя $(i)$  по параметру $a,$ получим $(ii).$ Полагая в $(ii)$ $b=1$ и дифференцируя по $a,$
получим $(iii).$  \qed

\begin{lemma}{\rm (см.~\cite[\S 2.8]{be})}  \label{l2}
Пусть $a, b, c\in {\mathbb R},$ $n\in {\mathbb N},$ $c$ отлично от отрицательных целых чисел и нуля,
и $c-a$ отлично от чисел $0, -1, -2, \ldots, -n+1.$
Тогда
$$
(i)\quad\frac{d^n}{dz^n}[z^{c-a+n-1}(1-z)^{a+b-c}
F(a,\, b;\, c;\, z)]=(c-a)_nz^{c-a-1}(1-z)^{a+b-c-n}F(a-n,\, b;\, c;\, z).
$$
В частности, при $a=0$ имеем
$$
\frac{d^n}{dz^n}[z^{n+c-1}(1-z)^{b-c}]=(c)_nz^{c-1}(1-z)^{b-c-n}F(-n,\, b;\, c;\, z).
$$
$(ii)$ Если и $c-b$ отлично от чисел $0, -1, -2, \ldots, -n+1,$ то
$$
\frac{d^n}{dz^n}[(1-z)^{a+b-c}
F(a,\, b;\, c;\, z)]=\frac{(c-a)_n(c-b)_n}{(c)_n}(1-z)^{a+b-c-n}F(a,\, b;\, c+n;\, z).
$$
\end{lemma}


\section{Рациональные приближения для $\ln(b)-\psi(a+1).$}


Пусть $a, b\in {\mathbb R},$ $b>0,$ $a>-1.$ Для целого $n\ge 0$ определим нтеграл
\begin{equation}
E_n=E_n^{(a,b)}=\frac{b^{2n+a+1}}{\Gamma(a+1)}\int_0^{\infty}Q_n^{(a,b)}(x)x^ae^{-bx}\ln(x)\, dx.
\label{eq04}
\end{equation}
\begin{lemma} \label{l3}
Для любого целого $n\ge 0$ справедливо следующее представление
$$
E_n=q_n(\psi(a+1)-\ln b)+p_n,
$$
где последовательности $p_n$ и $q_n$ определены в теореме {\rm\ref{t1}.}
\end{lemma}
{\bf Доказательство.} При $n=0$ имеем $Q_0^{(a,b)}\equiv 1,$ и
$E_0=\psi(a+1)-\ln(b).$ Пусть $n\ge 1.$
Из равенств (\ref{eq01}) и (\ref{eq04}) с помощью последовательного
ин\-тег\-ри\-ро\-ва\-ния по частям получаем
\begin{equation}
E_n=\frac{(-1)^n}{n!^2}\frac{b^{2n+a+1}}{\Gamma(a+1)}\int_0^{\infty}
\left(\frac{\ln(x)}{1-x}\right)^{(n)}x^n(x^{n+a}(1-x)^{2n+1}e^{-bx})^{(n)}\,dx.
\label{eq05}
\end{equation}
Вычисляя $n$-ую производную функции $\ln(x)/(1-x)$ по правилу Лейбница,
находим
$$
\left(\frac{\ln(x)}{1-x}\right)^{(n)}=\sum_{k=0}^n\binom{n}{k}(\ln(x))^{(k)}\left(\frac{1}{1-x}
\right)^{(n-k)}=\frac{n!\,\ln(x)}{(1-x)^{n+1}}-\sum_{k=1}^n\frac{(-1)^k n!}{kx^k(1-x)^{n-k+1}}.
$$
Подставляя найденное выражение в (\ref{eq05}) и интегрируя по частям, находим
\begin{equation}
\begin{split}
E_n&=\frac{b^{2n+a+1}}{\Gamma(a+1)n!}
\int_0^{\infty}\left[\left(\frac{x^n\ln x}{(1-x)^{n+1}}
\right)^{(n)}-\sum_{k=1}^n\frac{(-1)^k}{k}\left(\frac{x^{n-k}}{(1-x)^{n-k+1}}\right)^{(n)}\right]\\[5pt]
&\times x^{n+a}(1-x)^{2n+1}e^{-bx}\,dx.
\label{eq06}
\end{split}
\end{equation}
Применяя правило Лейбница и лемму \ref{l2}, вычислим $n$-ые производные в
формуле (\ref{eq06}). Имеем
\begin{equation*}
\begin{split}
&\left(\frac{x^n\ln(x)}{(1-x)^{n+1}}\right)^{(n)}=\sum_{k=0}^n\binom{n}{k}
(\ln(x))^{(k)}\left(\frac{x^n}{(1-x)^{n+1}}\right)^{(n-k)}\\
&=\ln(x) \left(
\frac{x^n}{(1-x)^{n+1}}\right)^{(n)}-\sum_{k=1}^n\binom{n}{k}
\frac{(-1)^k (k-1)!}{x^k}\left(\frac{x^n}{(1-x)^{n+1}}\right)^{(n-k)}\\
&=\frac{n! \ln(x)}{(1-x)^{2n+1}}F(-n,\, -n;\, 1;\, x)
-\sum_{k=1}^n\binom{n}{k}\frac{(-1)^k n!}{k(1-x)^{2n+1-k}}F(k-n,\, k-n;\, k+1;\, x).
\end{split}
\end{equation*}
Аналогично, по лемме \ref{l2} находим
\begin{equation*}
\begin{split}
\left(\frac{x^{n-k}}{(1-x)^{n-k+1}}\right)^{(n)}&=(n-k)!\left(
(1-x)^{2k-2n-1}F(k-n,\, k-n;\, 1;\, x)\right)^{(k)}\\
&=n!\binom{n}{k}(1-x)^{k-2n-1}F(k-n,\, k-n;\, k+1;\, x).
\end{split}
\end{equation*}
Тогда
\begin{equation}
\begin{split}
E_n&=\frac{b^{2n+a+1}}{\Gamma(a+1)}\left[\int_0^{\infty}x^{n+a}e^{-bx}\ln(x)
F(-n,\, -n;\, 1;\, x)\,dx\right.\\
&\left.-2\sum_{k=1}^n\frac{(-1)^k}{k}\binom{n}{k}\int_0^{\infty}
x^{n+a}(1-x)^ke^{-bx}F(k-n,\, k-n;\, k+1;\, x)\,dx\right].
\label{eq07}
\end{split}
\end{equation}
Интегралы в (\ref{eq07}) могут быть явно вычислены по лемме \ref{l1}. Действительно,
\begin{equation}
\begin{split}
 \int_0^{\infty}x^{n+a}e^{-bx}&\ln(x)
 F(-n,\, -n;\, 1;\, x)\,dx= \sum_{k=0}^n\binom{n}{k}^2
 \int_0^{\infty}x^{n+k+a}e^{-bx}\ln(x)\,dx \\
 &=\sum_{k=0}^n\binom{n}{k}^2\frac{\Gamma(n+k+a+1)}{b^{n+k+a+1}}
 (\psi(n+k+a+1)-\ln(b))
 \label{eq08}
 \end{split}
 \end{equation}
 Для второго интеграла имеем
 \begin{equation}
 \begin{split}
 &\int_0^{\infty}x^{n+a}(1-x)^ke^{-bx}F(k-n,\, k-n;\, k+1;\, x)\,dx\\
 &=\sum_{m=0}^k\sum_{l=0}^{n-k}(-1)^m\binom{k}{m}\frac{(k-n)_l^2}{l!(k+1)_l}\int_0^{\infty}
 x^{n+a+m+l}e^{-bx}\,dx\\
 &=\binom{n}{k}^{-1}\sum_{m=0}^k\sum_{l=0}^{n-k}(-1)^m\binom{k}{m}\binom{n-k}{l}\binom{n}{k+l}
 \frac{\Gamma(m+n+l+a+1)}{b^{m+n+l+a+1}}.
 \label{eq09}
 \end{split}
 \end{equation}
Наконец, из (\ref{eq07})--(\ref{eq09}) и функционального соотношения для дигамма-функции
$$
\psi(z+n)=\psi(z)+H_n(z-1)
$$
получаем требуемое утверждение.   \qed

\begin{lemma}   \label{l4}
Пусть $a,b \in{\mathbb R},$ $a>-1,$ $b>0.$ Тогда для любого  целого $n\ge 0$  имеем
\begin{equation}
E_n=\frac{b^{2n+a+1}}{\Gamma(a+1)}\int_0^{\infty} \int_0^{\infty}
\frac{x^{n+a}y^n(x-1)^{2n+1}e^{-bx}}{(xy+1)^{n+1}(y+1)^{n+1}}\,dxdy.
\label{En}
\end{equation}
\end{lemma}

 {\bf Доказательство.}  Подставим выражение (\ref{eq01})  в (\ref{eq04}), получим
 $$
 E_n=\frac{b^{2n+a+1}}{\Gamma(a+1) n!^2}\int_0^{\infty}\frac{\ln(x)}{1-x}(x^n(x^{n+a}(1-x)^{2n+1}
 e^{-bx})^{(n)})^{(n)}\,dx.
 $$
 Воспользуемся представлением
 $$
 \frac{\ln(x)}{1-x}=-\int_0^{\infty}\frac{du}{(1+u)(1+ux)}
 $$
и осуществим  $n$-кратное интегрирование по частям, получим
$$
E_n=-\frac{b^{2n+a+1}}{\Gamma(a+1) n!}\int_0^{\infty}\int_0^{\infty}\frac{x^nu^n}{(1+u)(1+xu)^{n+1}}
(x^{n+a}(1-x)^{2n+1}e^{-bx})^{(n)}\,dxdu.
$$
Сделаем замену переменной $u$ по формуле $u=1/(xy)$ или $y=1/(xu),$ имеем
$$
E_n=-\frac{b^{2n+a+1}}{\Gamma(a+1) n!}\int_0^{\infty}\int_0^{\infty}
\frac{(x^{n+a}(1-x)^{2n+1}e^{-bx})^{(n)}}{(xy+1)(y+1)^{n+1}}\,dxdy.
$$
Снова интегрируя по частям, получим
$$
E_n=-\frac{b^{2n+a+1}}{\Gamma(a+1)}\int_0^{\infty}\int_0^{\infty}
\frac{x^{n+a}y^n(1-x)^{2n+1}e^{-bx}}{(y+1)^{n+1}(xy+1)^{n+1}}\, dxdy,
$$
и лемма доказана. \qed

{\it Замечание.} Применяя к последовательности интегралов $E_n^{(a,b)}$ алгоритм
MultiAlmkvistZeilberger (теоретическое обоснование которого можно найти в \cite{maz}),
не\-труд\-но получить, что последовательность $E_n^{(a,b)}$ удовлетворяет рекуррентному
урав\-не\-нию третьего порядка с полиномиальными коэффициентами. С помощью ал\-го\-рит\-ма EKHAD
можно показать, что последовательность $q_n=\sum_{k=0}^n\binom{n}{k}^2(a+1)_{n+k}b^{n-k}$
удовлетворяет тому же рекуррентному уравнению, и следовательно, по линейности последовательность
$p_n$ также является решением того же уравнения. Мы не вы\-пи\-сы\-ва\-ем это уравнение в общем виде
в силу громоздкости его коэффициентов. Однако, при конкретных значениях $a$ и $b$ это легко можно сделать
(см.,например, замечания к следствиям \ref{c1}, \ref{c3}, \ref{c4}).
\begin{lemma} \label{l5}
Справедливы асимптотические формулы
\begin{equation*}
\begin{split}
q_n&=(2n)!\frac{e^{\sqrt{2bn}}}{n^{1/4-a}}(c(a,b)+O(n^{-1/2})), \qquad\quad\, n\to\infty, \\
E_n&=(2n)!\frac{e^{-\sqrt{2bn}}}{n^{1/4-a}}(2\pi c(a,b)+O(n^{-1/2})), \qquad n\to\infty,
\end{split}
\end{equation*}
где $c(a,b)=2^a/(\sqrt{\pi}\,\Gamma(a+1)e^{3b/8}(8b)^{1/4}).$
\end{lemma}

{\bf Доказательство.} Асимптотику последовательности $q_n$ найдем, следуя методу
работы \cite{apl}. Используя лемму \ref{l1} и проведя последовательное интегрирование
по частям, находим
\begin{equation*}
\begin{split}
q_n&=\frac{b^{2n+a+1}}{\Gamma(a+1)}\int_0^{\infty}Q_n^{(a,b)}(x)x^ae^{-bx}\,dx \\
&=\frac{(-1)^n}{n!}\frac{b^{2n+a+1}}{\Gamma(a+1)}\int_0^{\infty}x^{n+a}(x-1)^{2n+1}
e^{-bx}\left(\frac{x^n}{(x-1)^{n+1}}\right)^{(n)}\,dx.
\end{split}
\end{equation*}
Далее дословно повторяя рассуждения из \cite[стр.~56-59]{apl}, получим
\begin{equation}
q_n=\frac{b^{2n+a+1}n^{2n+a+3/4}}{2\sqrt{\pi}\,\Gamma(a+1)}\int_0^{\infty}
e^{n\varphi_1(y)}\psi_1(y)\,dy\, (1+O(n^{-1/2}),
\label{eq10}
\end{equation}
где
$$
\varphi_1(y)=2\ln(\sqrt{y}+1/\sqrt{n})+\ln y- by, \qquad
\psi_1(y)=(\sqrt{y}+1/\sqrt{n})y^{a-1/4}.
$$
Функция $\varphi_1(y)$ достигает глобального максимума на полуоси $(0,+\infty)$
в точке $y_0=2/b-1/\sqrt{2bn}+3/(8n)+O(n^{-3/2}).$  При этом имеем
\begin{equation*}
\begin{split}
\varphi_1(y_0)&=2\ln(2/b)-2+\sqrt{\frac{2b}{n}}-\frac{3b}{8n}+O(n^{-3/2}), \\
\varphi''_1(y_0)&=-\frac{b^2}{2}+O(n^{-1/2}), \\
\psi_1(y_0)&=\left(\frac{2}{b}\right)^{a+1/4}+O(n^{-1/2}).
\end{split}
\end{equation*}
Применяя метод перевала к интегралу (\ref{eq10}), получим асимптотическую формулу 
$$
q_n=\frac{(2n)^{2n+a+1/4}}{\Gamma(a+1)b^{1/4}}\,
e^{-2n+\sqrt{2bn}-3b/8}(1+O(n^{-1/2}))=(2n)!\frac{e^{\sqrt{2bn}}}{n^{1/4-a}}(c(a,b)+O(n^{-1/2})).
$$
Асимптотику остатков приближений $E_n$ найдем, используя лемму \ref{l4} и метод Лапласа
для кратных интегралов (см.~\cite[\S 4]{fe}). Сделаем в интеграле (\ref{En}) замену переменных
$x=nu,$ $y=v/\sqrt{nu},$ получим
\begin{equation}
E_n=\frac{b^{2n+a+1}n^{3n+a+2}}{\Gamma(a+1)}\int_0^{\infty}\int_0^{\infty}
g(u,v) f^n(u,v)\,dudv,
\label{eq11}
\end{equation}
где
$$
f(u,v)=\frac{uv(u-1/n)^2e^{-bu}}{(v+\sqrt{nu})(v\sqrt{nu}+1)}, \qquad
g(u,v)=\frac{u^a(u-1/n)}{(v+\sqrt{nu})(v\sqrt{nu}+1)}.
$$
Как нетрудно видеть, функция $f(u,v)$ достигает глобального максимума в области
$u>0,$ $v>0$ в точке $(u_0, v_0),$ где $u_0=2/b+1/\sqrt{2bn}+3/(8n)+O(n^{-3/2}),$ $v_0=1.$
При этом имеем
$$
f(u_0,v_0)=\frac{u_0}{n^2}(\sqrt{nu_0}-1)^2e^{-bu_0}=
\exp(2\ln(2/b)-\ln n-2-\sqrt{2b/n}-3b/(8n)+O(n^{-3/2})),
$$
$$
\frac{\partial^2f}{\partial u^2}(u_0,v_0)=-\frac{2e^{-2}}{n}(1+O(n^{-1/2})),
\qquad
\frac{\partial^2f}{\partial v^2}(u_0,v_0)=-\frac{4\sqrt{2}\,e^{-2}}{(bn)^{3/2}}(1+O(n^{-1/2})),
$$
$$
\frac{\partial^2f}{\partial u \partial v}(u_0,v_0)=
\frac{\partial^2f}{\partial v \partial u}(u_0,v_0)=0, \quad
g(u_0,v_0)=\frac{1}{n}\left(\frac{2}{b}\right)^a(1+O(n^{-1/2})).
$$
Применяя метод Лапласа к интегралу (\ref{eq11}), получим
$$
E_n=\frac{b^{2n+a+1}n^{3n+a+2}}{\Gamma(a+1)}\frac{2\pi}{n\sqrt{A}}\, g(u_0,v_0)
f^{n+1}(u_0,v_0)(1+O(n^{-1})),
$$
где $A=\det\left(\begin{smallmatrix}
\frac{\partial^2f}{\partial u^2}(u_0,v_0)&\frac{\partial^2f}{\partial u \partial v}(u_0,v_0)\\
\frac{\partial^2f}{\partial v \partial u}(u_0,v_0)&\frac{\partial^2f}{\partial v^2}(u_0,v_0)
\end{smallmatrix} \right),$
и следовательно,
$$
E_n=\frac{2\pi (2n)^{2n+a+1/4}}{\Gamma(a+1) b^{1/4}}e^{-2n-\sqrt{2bn}-3b/8}(1+O(n^{-1/2}))=
(2n)!\frac{e^{-\sqrt{2bn}}}{n^{1/4-a}}(2\pi c(a,b)+O(n^{-1/2})).
$$
Лемма доказана. \qed

Теорема \ref{t1} легко следует из лемм \ref{l3}--\ref{l5}.

\section{Свойства обобщенных многочленов Лагерра.}

Пусть $a_1,\ldots, a_s\in {\mathbb R},$
$a_1, \ldots, a_s>-1$ и $n_1,\ldots, n_s$ целые неотрицательные числа. Определим
обобщенные многочлены Лагерра $L_s(x)=L_{(n_1,\ldots,n_s)}^{(a_1,\ldots,a_s)}$ рекурсивно.
По\-ло\-жим $L_0(x)\equiv 1$ и
\begin{equation}
L_s(x)=\frac{x^{-a_s}e^x}{n_s!}(x^{n_s+a_s}e^{-x}L_{s-1}(x))^{(n_s)}, \quad s\ge 1.
\label{eq101}
\end{equation}
\begin{lemma}  \label{l6}
Для обобщенных многочленов Лагерра имеет место  представление
$$
L_{s}(x)=\sum_{k_1=0}^{n_1}\ldots\sum_{k_s=0}^{n_s}\prod_{j=1}^s
\binom{n_j+a_j+k_1+\cdots+k_{j-1}}{n_j-k_j}\frac{(-x)^{k_j}}{k_j!}.
$$
\end{lemma}
{\bf Доказательство.}  Доказательство проведем индукцией по $s.$ При $s=0$ имеем $L_{0}(x)\equiv 1.$
При $s=1$ получаем классические многочлены Лагерра
$$ L_{1}(x)=L_{n_1}^{a_1}(x)=\sum_{k=0}^{n_1}\binom{n_1+a_1}{n_1-k}\frac{(-x)^k}{k!}.
$$
Пусть $s>1,$ и предположим, что  утверждение  верно до $s-1$ включительно.
Тогда имеем
\begin{equation*}
\begin{split}
 L_{s}(x)&=\frac{x^{-a_s}e^x}{n_s!}(x^{n_s+a_s}e^{-x}L_{s-1}(x))^{(n_s)}
=\frac{x^{-a_s}e^x}{n_s!}\\
&\times\left(\sum_{k_1=0}^{n_1}\ldots\sum_{k_{s-1}=0}^{n_{s-1}}\prod_{j=1}^{s-1}
\frac{(-1)^{k_j}}{k_j!}\binom{n_j+a_j+k_1+\cdots+k_{j-1}}{n_j-k_j} x^{n_s+a_s+k_1+\cdots
+k_{s-1}}e^{-x}\right)^{(n_s)}\\
&=\frac{1}{n_s!}\sum_{k_1=0}^{n_1}\ldots\sum_{k_{s-1}=0}^{n_{s-1}}\sum_{k_s=0}^{n_s}
\binom{n_s}{k_s}
\prod_{j=1}^{s-1}
\frac{(-x)^{k_j}}{k_j!}\binom{n_j+a_j+k_1+\cdots+k_{j-1}}{n_j-k_j}\times (-x)^{k_s}\\[3pt]
&\times\frac{(n_s+a_s+\sum_{j=1}^{s-1}k_j)!}{(k_1+\cdots+k_s)!}
=\sum_{k_1=0}^{n_1}\ldots\sum_{k_s=0}^{n_s}\prod_{j=1}^s
\binom{n_j+a_j+k_1+\cdots+k_{j-1}}{n_j-k_j}\frac{(-x)^{k_j}}{k_j!},
\end{split}
 \end{equation*}
и лемма доказана.  \qed
\begin{lemma} \label{l6.5}
Пусть $s, n_1, \ldots, n_s\in {\mathbb N}, a_1, \ldots, a_s\in {\mathbb R},$
$a_1, \ldots, a_s>-1,$ и $a_j-a_k\not\in {\mathbb Z}$ при $j\ne k.$
Для обобщенных многочленов Лагерра имеют место следующие соотношения ортогональности
\begin{equation}
\int_0^{\infty}t^{l+a_m}L_s(t)e^{-t}\,dt=0, \quad l=0,1,\ldots, n_m-1,\quad m=1,\ldots,s.
\label{eq102}
\end{equation}
Кроме того, при целом $l\ge n_m$ имеем
\begin{equation}
\int_0^{\infty}t^{l+a_m}L_s(t)e^{-t}\,dt=(-1)^{n_1+\ldots+n_s}\Gamma(l+a_m+1)\prod_{j=1}^s
\binom{l+a_m-a_j}{n_j}.
\label{eq103}
\end{equation}
\end{lemma}
{\bf Доказательство.} Как легко видеть из (\ref{eq101}), для многочленов $L_s(x)$
справедлива формула Родрига
$$
L_s(x)=e^x\prod_{j=1}^s\left[\frac{x^{-a_j}}{n_j!}\frac{d^{n_j}}{dx^{n_j}}x^{n_j+a_j}\right]e^{-x},
$$
где операторы, стоящие в квадратных скобках, коммутируют между собой (см.~\cite[\S 2.3]{as}).
Используя свойство коммутирования и $n_m$-кратное интегрирование по час\-тям, легко получаем (\ref{eq102}).
Равенство (\ref{eq103}) непосредственно следует из формулы ин\-тег\-ри\-ро\-ва\-ния по частям. \qed

Для многочленов $A_{s,n}(x)=L_{(n,\ldots,n)}^{(0,\ldots,0)}(x),$ определенных в
формулировке теоремы \ref{t2}, имеем.
\begin{lemma}\cite[\S 4.2]{ri}   \label{l7}
Пусть $s, n\in {\mathbb N}.$ Для обобщенных многочленов Лагерра $A_{s,n}$
 имеют место следующие соотношения ортогональности
\begin{equation}
\int_0^{\infty}t^lA_{s,n}(t)\ln(t)^{j-1}e^{-t}\,dt=0, \qquad l=0,1,\ldots, n-1,
\quad j=1,2,\ldots,s.
\label{eq12}
\end{equation}
\end{lemma}
{\bf Доказательство.} Легко следует из формулы интегрирования по частям.  \qed
\begin{lemma} \label{l8}
При $s,l,n\in {\mathbb N}\cup\{0\}$   и $l\ge n$  имеем
\begin{equation*}
\begin{array}{lrcl}
(a)&\displaystyle\frac{(-1)^{sn}}{l!}\int_0^{\infty}t^lA_{s,n}(t)e^{-t}\,dt\!\!&\!=\!&\!\!\displaystyle
 \binom{l}{n}^s,    \\[12pt]
(b)&\displaystyle\frac{(-1)^{sn}}{l!}\int_0^{\infty}t^lA_{s,n}(t)\ln(t)e^{-t}\,dt\!\!&\!=\!&\!\!\displaystyle
 \binom{l}{n}^s
\psi(l+1)-s\binom{l}{n}^{s-1}\sum_{k=1}^n\frac{(-1)^k}{k}\binom{l}{n-k}, \\[15pt]
(c)&\displaystyle\frac{(-1)^{sn}}{l!}\int_0^{\infty}t^lA_{s,n}(t)\ln^2(t)e^{-t}\,dt\!\!&\!=\!&\!\!\displaystyle
\binom{l}{n}^s(\psi^2(l+1)+\psi'(l+1))-2s\binom{l}{n}^{s-1}
\end{array}
\end{equation*}
$$
\times\sum_{k=1}^n
\frac{(-1)^k}{k}\binom{l}{n-k}(\psi(l+1)+H_{k-1})+ s(s-1)\binom{l}{n}^{s-2}
\left(\sum_{k=1}^n\frac{(-1)^k}{k}\binom{l}{n-k}\right)^2.
$$
\end{lemma}
{\bf Доказательство.} При $s=0$  многочлен $A_{0,n}(t)$  тождественно равен $1,$
и тре\-буе\-мые утверждения легко следуют из леммы \ref{l1}. При $s\ge 1$ по формуле интегрирования
по частям имеем
\begin{equation} I_{s,j}:=(-1)^{sn}\int_0^{\infty}t^lA_{s,n}(t)\ln(t)^{j-1}e^{-t}dt
=\frac{(-1)^{(s-1)n}}{n!}\int_0^{\infty}(t^l\ln(t)^{j-1})^{(n)}t^nA_{s-1,n}(t)e^{-t}dt.
\label{eq13}
\end{equation}
Откуда, в частности, получаем $I_{s,1}=\binom{l}{n}I_{s-1,1},$ и формула $(a)$ доказана. Для
доказательства $(b)$ вычислим $n$-ую производную произведения $t^l\ln(t)$ по правилу Лейбница,
имеем
$$
(t^l\ln(t))^{(n)}=\sum_{k=0}^n\binom{n}{k}(t^l)^{(n-k)}(\ln(t))^{(k)}=
n! t^{l-n}\left(\binom{l}{n}\ln(t)-\sum_{k=1}^n\frac{(-1)^k}{k}\binom{l}{n-k}\right).
$$
И следовательно, с учетом (\ref{eq13}) и  $(a)$  находим
$$
I_{s,2}=\binom{l}{n}I_{s-1,2}-I_{s-1,1}\cdot \sum_{k=1}^n\frac{(-1)^k}{k}\binom{l}{n-k},
$$
откуда формула $(b)$  легко следует индукцией по $s.$
Для доказательства $(c),$ ана\-ло\-гич\-но, по правилу Лейбница находим
$$
(t^l\ln(t)^2)^{(n)}=n! t^{l-n}\left(\binom{l}{n}\ln(t)^2-2\sum_{k=1}^n
\frac{(-1)^k}{k}\binom{l}{n-k}(\ln(t)+H_{k-1})\right).
$$
Тогда из (\ref{eq13}) получаем
$$
I_{s,3}=\binom{l}{n}I_{s-1,3}-2\sum_{k=1}^n\frac{(-1)^k}{k}\binom{l}{n-k}
(I_{s-1,2}+H_{k-1}I_{s-1,1}),
$$
и требуемое утверждение легко следует по индукции по параметру $s.$  \qed

Из леммы \ref{l8} легко получаем следующее утверждение.
\begin{lemma} \label{l9}
При $s,l,n\in {\mathbb N}\cup\{0\}$   и $l\ge n$
$$
(-1)^{sn}\int_0^{\infty}t^lA_{s,n}(t)(\ln(t)^2+2\gamma\ln(t))e^{-t}\,dt=
p_{l,n}(\zeta(2)-\gamma^2)+q_{l,n},
$$
где $p_{l,n}=l! \binom{l}{n}^s$ и $\frac{D_n^2}{n!}q_{l,n}\in {\mathbb Z}.$
\end{lemma}

\section{Доказательство теоремы 2.}

Разложим интеграл $J_{s,n}(x)$ в линейную форму от чисел $1$  и $\gamma,$
используя лемму \ref{l8}. Для этого запишем $J_{s,n}$ в виде
\begin{equation}
\begin{split}
\frac{J_{s,n}(x)}{(n+1)^s}
&=\int_0^{\infty}A_{s,n+1}(t)\frac{A_{s,n}(t)-A_{s,n}(x)}{t-x}
\ln(t)e^{-t}\,dt \\
&-\int_0^{\infty}A_{s,n}(t)\frac{A_{s,n+1}(t)-A_{s,n+1}(x)}{t-x}
\ln(t)e^{-t}\,dt.
\label{eq14}
\end{split}
\end{equation}
По лемме  \ref{l6} преобразуем частное двух многочленов, имеем
$$
\frac{A_{s,n}(t)-A_{s,n}(x)}{t-x}=\sum_{k_1=0}^n\ldots\sum_{k_s=0}^n
\prod_{j=1}^s\frac{(-1)^{k_j}}{k_j!}\binom{n+k_1+\ldots+k_{j-1}}{n-k_j}
\!\sum_{l=0}^{k_1+\ldots+k_s-1}\!\!\!t^lx^{k_1+\ldots+k_s-l-1}.
$$
Используя соотношения ортогональности (\ref{eq12}) и лемму \ref{l8}, получаем
\begin{equation}
\begin{split}
&\int_0^{\infty}A_{s,n+1}(t)\frac{A_{s,n}(t)-A_{s,n}(x)}{t-x}\ln(t)e^{-t}\,dt
=\underset{k_1+\ldots+k_s\ge n+2}{\sum_{k_1=0}^n\ldots\sum_{k_s=0}^n}
\sum_{l=n+1}^{k_1+\ldots+k_s-1}\prod_{j=1}^s
\frac{(-1)^{k_j}}{k_j!} \\
&\times\binom{n+k_1+\ldots+k_{j-1}}{n-k_j}x^{k_1+\ldots+k_s-l-1}
\int_0^{\infty}t^lA_{s,n+1}(t)\ln(t)e^{-t}dt \\
&=\pm\underset{k_1+\ldots+k_s\ge n+2}{\sum_{k_1=0}^n\ldots\sum_{k_s=0}^n}
\sum_{l=n+1}^{k_1+\ldots+k_s-1}\prod_{j=1}^s
\frac{(-1)^{k_j}}{k_j!}\binom{n+k_1+\ldots+k_{j-1}}{n-k_j}x^{k_1+\ldots+k_s-l-1}
l! \\
&\times\left(\binom{l}{n+1}^s(H_l-\gamma)-s\binom{l}{n+1}^{s-1}\sum_{k=1}^{n+1}
\frac{(-1)^k}{k}\binom{l}{n+1-k}\right) \\
&\in
\frac{1}{n!^{s-1}}\Bigl(\frac{1}{D_n}{\mathbb Z}[x]+\gamma{\mathbb Z}[x]\Bigr).
\label{eq15}
\end{split}
\end{equation}
Аналогично, для второго интеграла имеем
\begin{equation}
\begin{split}
&\int_0^{\infty}A_{s,n}(t)\frac{A_{s,n+1}(t)-A_{s,n+1}(x)}{t-x}\ln(t)e^{-t}\,dt\\
&=\pm\underset{k_1+\ldots+k_s\ge n+1}{\sum_{k_1=0}^{n+1}\ldots\sum_{k_s=0}^{n+1}}\sum_{l=n}^{k_1+\ldots+k_s-1}\prod_{j=1}^s
\frac{(-1)^{k_j}}{k_j!}\binom{n+1+k_1+\ldots+k_{j-1}}{n+1-k_j}x^{k_1+\ldots+k_s-l-1}
l! \\
&\times\left(\binom{l}{n}^s(H_l-\gamma)-s\binom{l}{n}^{s-1}\sum_{k=1}^{n}
\frac{(-1)^k}{k}\binom{l}{n-k}\right)
\in
\frac{1}{(n+1)^sn!^{s-1}}\Bigl(\frac{1}{D_n}{\mathbb Z}[x]+\gamma{\mathbb Z}[x]\Bigr).
\label{eq16}
\end{split}
\end{equation}
Наконец, из (\ref{eq14})--(\ref{eq16}) получаем
$
J_{s,n}(x)=Q_{s,n}(x)\gamma-P_{s,n}(x),
$
где $n!^{s-1}Q_{s,n}(x)\in {\mathbb Z}[x],$
$D_nn!^{s-1}P_{s,n}(x)\in {\mathbb Z}[x],$ и теорема доказана.  \qed


При $s=2$ из теоремы \ref{t2} получаем следствие \ref{c5}.
При $s=3$ получаем  первых два включения следствия \ref{c6}.
Для доказательства третьего включения разложим интеграл, стоящий в левой части
формулы (\ref{eqri2}),  по лемме \ref{l9}. Имеем
\begin{equation*}
\begin{split}
&\int_0^{\infty}\frac{A_{3,n}(t)A_{3,n+1}(1)-A_{3,n+1}(t)A_{3,n}(1)}{1-t}
(\ln(t)^2+2\gamma\ln(t))e^{-t}dt \\
&=\int_0^{\infty}A_{3,n}(t)\frac{A_{3,n+1}(1)-A_{3,n+1}(t)}{1-t}
(\ln(t)^2+2\gamma\ln(t))e^{-t}dt \\
&-\int_0^{\infty}A_{3,n+1}(t)\frac{A_{3,n}(1)-A_{3,n}(t)}{1-t}
(\ln(t)^2+2\gamma\ln(t))e^{-t}dt \\
&=\underset{k_1+k_2+k_3\ge n+1}{\sum_{k_1=0}^{n+1}\sum_{k_2=0}^{n+1}%
\sum_{k_3=0}^{n+1}}\prod_{j=1}^3\frac{(-1)^{k_j}}{k_j!}\binom{n+1+k_1+\ldots+k_{j-1}}{n+1-k_j}
\sum_{l=n}^{k_1+k_2+k_3-1}\int_0^{\infty}t^lA_{3,n}(t) \\
&\times(\ln(t)^2+2\gamma\ln(t))e^{-t}dt-
\underset{k_1+k_2+k_3\ge n+2}{\sum_{k_1=0}^{n}\sum_{k_2=0}^{n}%
\sum_{k_3=0}^{n}}\prod_{j=1}^3\frac{(-1)^{k_j}}{k_j!}\binom{n+k_1+\ldots+k_{j-1}}{n-k_j} \\
&\times\sum_{l=n+1}^{k_1+k_2+k_3-1}\int_0^{\infty}t^lA_{3,n+1}(t)
(\ln(t)^2+2\gamma\ln(t))e^{-t}dt
=\underset{k_1+k_2+k_3\ge n+1}{\sum_{k_1=0}^{n+1}\sum_{k_2=0}^{n+1}%
\sum_{k_3=0}^{n+1}}\prod_{j=1}^3\frac{(-1)^{k_j+n}}{k_j!} \\
&\times\binom{n+1+k_1+\ldots+k_{j-1}}{n+1-k_j}
\sum_{l=n}^{k_1+k_2+k_3-1}(p_{l,n}(\zeta(2)-\gamma^2)+q_{l,n}) \\
&-\underset{k_1+k_2+k_3\ge n+2}{\sum_{k_1=0}^{n}\sum_{k_2=0}^{n}%
\sum_{k_3=0}^{n}}\prod_{j=1}^3\frac{(-1)^{k_j+n}}{k_j!}\binom{n+k_1+\ldots+k_{j-1}}{n-k_j}
\sum_{l=n+1}^{k_1+k_2+k_3-1}(p_{l,n+1}(\zeta(2)-\gamma^2)+q_{l,n+1}) \\
&\in \frac{1}{(n+1)^3n!^2}\Bigl(\frac{1}{D_n^2}{\mathbb Z}+(\zeta(2)-\gamma^2){\mathbb Z}\Bigr),
\end{split}
\end{equation*}
и следствие \ref{c6} полностью доказано. \qed

\section{Доказательство предложения 1.}

Из равенства (\ref{eq100}) легко следует, что
$$
Q_n(a,x)=\frac{1}{\Gamma(a+1)}\int_0^{+\infty}\frac{\hat{A}_{n,a}(x)\hat{A}_{n+1,a}(t)
-\hat{A}_{n+1,a}(x)\hat{A}_{n,a}(t)}{t-x}t^ae^{-t}\,dt
$$
и
$$
P_n(a,x)=\int_0^{+\infty}\frac{\hat{A}_{n,a}(x)\hat{A}_{n+1,a}(t)
-\hat{A}_{n+1,a}(x)\hat{A}_{n,a}(t)}{t-x}e^{-t}\,dt.
$$
Перепишем  многочлен $Q_n(a,x)$ в виде
\begin{equation*}
\begin{split}
Q_n(a,x)&=\frac{1}{\Gamma(a+1)}\int_0^{\infty}\hat{A}_{n,a}(t)
\frac{\hat{A}_{n+1,a}(x)
-\hat{A}_{n+1,a}(t)}{x-t}t^ae^{-t}\,dt \\
&-\frac{1}{\Gamma(a+1)}\int_0^{\infty}\hat{A}_{n+1,a}(t)
\frac{\hat{A}_{n,a}(x)-\hat{A}_{n,a}(t)}{x-t}t^ae^{-t}\,dt=J_2-J_1.
\end{split}
\end{equation*}
Так как
$$
\hat{A}_{n,a}(t)=L_{(n,n)}^{(a,0)}(t)=\sum_{k=0}^n\sum_{j=0}^n\binom{n+a}{n-k}
\binom{n+k}{n-j}\frac{(-x)^{k+j}}{k!j!},
$$
то применяя к интегралу $J_1$ лемму \ref{l6.5}, имеем
$$
J_1
=
\underset{k+j\ge n+2}{\sum_{k=0}^n\sum_{j=0}^n}\sum_{l=n+1}^{k+j-1}
\binom{n+a}{n-k}\binom{n+k}{n-j}\binom{l}{n+1}\binom{l+a}{n+1}\frac{(-1)^{k+j}}{k!j!}
(a+1)_lx^{k+j-1-l}.
$$
Следовательно, $J_1$ является многочленом от $x$ степени, не превосходящей $n-2.$ Более того,
используя хорошо известное свойство биномиальных коэффициентов (см., например, \cite[лемма 4.1]{chu})
$$
\mu_{a}^n \cdot \frac{(a)_n}{n!}\in {\mathbb Z},
$$
получаем $n!\,({\rm den}\, x)^{n-2}\mu_{a}^{3n}\, J_1\in {\mathbb Z}.$
Аналогично, для интеграла $J_2$ имеем
$$
J_2=\underset{k+j\ge n+1}{\sum_{k=0}^{n+1}\sum_{j=0}^{n+1}}\sum_{l=n}^{k+j-1}
\binom{n+1+a}{n+1-k}\binom{n+1+k}{n+1-j}\binom{l}{n}\binom{l+a}{n}\frac{(-1)^{k+j}}{k!\,j!}
(a+1)_lx^{k+j-l-1}.
$$
Откуда следует, что $(n+1)\cdot (n+1)!\,({\rm den}\, x)^{n+1}\mu_{a}^{3n+1}\,
J_2\in {\mathbb Z},$ и  требуемое включение для многочленов $Q_n(a,x)$ доказано.

Аналогичным образом, раскладывая интеграл, определяющий многочлен $P_n(a,x),$ имеем
$P_n(a,x)=J_4-J_3,$ где
$$
J_3=\underset{k+j\ge n+2}{\sum_{k=0}^n\sum_{j=0}^n}\sum_{l=n+1}^{k+j-1}
\binom{n+a}{n-k}\binom{n+k}{n-j}\binom{l-a}{n+1}\binom{l}{n+1}\frac{(-1)^{k+j}\,l!}{k!\,j!}x^{k+j-l-1},
$$
$$
J_4=\underset{k+j\ge n+1}{\sum_{k=0}^{n+1}\sum_{j=0}^{n+1}}\sum_{l=n}^{k+j-1}
\binom{n+1+a}{n+1-k}\binom{n+1+k}{n+1-j}\binom{l-a}{n}\binom{l}{n}\frac{(-1)^{k+j}\,l!}{k!\,j!}x^{k+j-l-1}.
$$
Откуда следует, что
$$
n!\, ({\rm den}\, x)^{n-2} \mu_{a}^{2n+1}\, J_3
\in {\mathbb Z}, \quad
(n+1)^2 n!\, ({\rm den}\, x)^{n+1} \mu_{a}^{2n+1}\, J_4
\in {\mathbb Z},
$$
и предложение доказано. \qed

\section{Заключительные замечания.}

В заключение, сравнивая обе конструкции Ап\-те\-ка\-ре\-ва и Ривоаля,
задающие   ра\-цио\-наль\-ные приближения
 константы Эйлера
\begin{equation}
\frac{p_n}{q_n}-\gamma=O(e^{-2\sqrt{2n}})
\quad\text{или}\quad \frac{D_n}{n!}p_n-\gamma\frac{D_n}{n!}q_n=O(4^nn^{n-1/4}e^{-\sqrt{2n}}),
\label{a1}
\end{equation}
и
\begin{equation}
\left|\gamma-\frac{P_n}{Q_n}\right|=O(e^{-9/2n^{2/3}+3/2n^{1/3}})\,\,\,\text{или}
\,\,\, |\gamma D_nn!\,Q_n-D_nn!\,P_n|=O(n^{n+1/2}e^{-3/2n^{2/3}+1/2n^{1/3}}),
\label{r1}
\end{equation}
где $p_n,$ $q_n$ определены в (\ref{znam}), а $P_n=P_n(1),$ $Q_n=Q_n(1)$ (см.~(\ref{pq})),
заметим, что конструкция Ривоаля (\ref{r1}) несколько ``ближе'' к диофантовым при\-бли\-же\-ни\-ям
в том смысле, что линейные формы в (\ref{a1}) и (\ref{r1}) имеют целые ко\-эф\-фи\-ци\-ен\-ты, и формы в (\ref{r1})
растут медленнее, чем формы в (\ref{a1}). С другой стороны, знаменатели в (\ref{a1}) выражаются в достаточно
компактной форме $\sum_{k=0}^n\binom{n}{k}^2(n+k)!$ и по своему виду напоминают числа Апери
$$
A_n=\sum_{k=0}^n\binom{n}{k}^2\binom{n+k}{k}, \qquad
B_n=\sum_{k=0}^n\binom{n}{k}^2\binom{n+k}{k}^2,
$$
которые являются знаменателями диофантовых приближений чисел $\zeta(2)$ и $\zeta(3)$
(см. \cite{po, pre}), позволяющих доказать  иррациональность $\zeta(2)$
и $\zeta(3).$  Чис\-ла Апери
удовлетворяют рекуррентным уравнениям второго порядка с по\-ли\-но\-ми\-аль\-ны\-ми коэффициентами.
Поэтому было бы интересно, с одной стороны, найти ра\-цио\-наль\-ные приближения константы
Эйлера, задаваемые рекуррентным уравнением второго порядка, с другой --- попытаться
улучшить скорость сходимости ра\-цио\-наль\-ных приближений, т.~е.~найти последователность
рациональных дробей $\tilde{p}_n/\tilde{q}_n,$ сходящихся к $\gamma$ геометрически.

\end{document}